\documentclass[12pt]{amsart}                          

\usepackage{tikz}
\usetikzlibrary{arrows,patterns,matrix,spy}

\begin{document}

\title{Progressions
 in reasoning in K-12 mathematics}
\author{Dev Sinha}
\maketitle 

``I think that it is  impossible for some of our students to learn to do proofs,'' explained a colleague of mine.  
My belated response is that all students can indeed learn to do proofs.  
College faculty have been asking aspiring math majors 
to make a huge jump, from not being responsible for providing
reasoning in entry-level college courses and most of their K-12 
experience to making formal proofs.  
To ease this transition, our department at the University of Oregon has recently created ``lab'' courses for
first-year students  in addition to our ``bridge'' requirement,
 to help students by degrees gain experience with proof-based mathematics.

More broadly, we should afford all students, at all levels, practice with gradually more demanding reasoning,
something which educators would call engaging in progressions in reasoning \cite{CCWT18}.  
Teachers, mathematics educators, and mathematicians
 have been working together to develop and study such learning progressions for 
 reasoning\footnote{We use the terms ``reasoning'' and ``proof'' interchangeably. There
is a continuum of formalism in reasoning and proof, and it is pedagogically useful
 for students to engage at  different levels, even within the same 
classroom activity.  We prefer to emphasize this range rather 
than making only a binary distinction based on some hard-to-define cutoff for rigor and formalism.}   in K-12 mathematics 
 (see for example \cite{SBK10, Knu02, KCSS02}).  
 In this article I share perspective 
 as a research mathematician who has developed and implemented reasoning-focused tasks
 for K-12 students of all grades, for both aspiring and current teachers, and for many types of undergraduates.\\

Historically, explicit calls for student reasoning had been all but absent in school curricula
until Euclidean geometry.  But mathematical
argument does not need to wait until the onset of puberty.  Indeed, children naturally ask ``why questions,'' and
based on research
both the National Council of Teachers of Mathematics Principles and Standards
for School Mathematics \cite{NCTM00} and Common Core State 
Standards in Mathematics,  in particular its third Practice Standard \cite{CCWT10},
call for more reasoning and proof throughout school mathematics.  

Such calls are being answered.  For example using
pictures such as the one below, many current curricula provide opportunities to reason about the sum of two odd numbers 
in the second or third grade.  
 
 \begin{center}\includegraphics[keepaspectratio=true, width=0.5\textwidth]{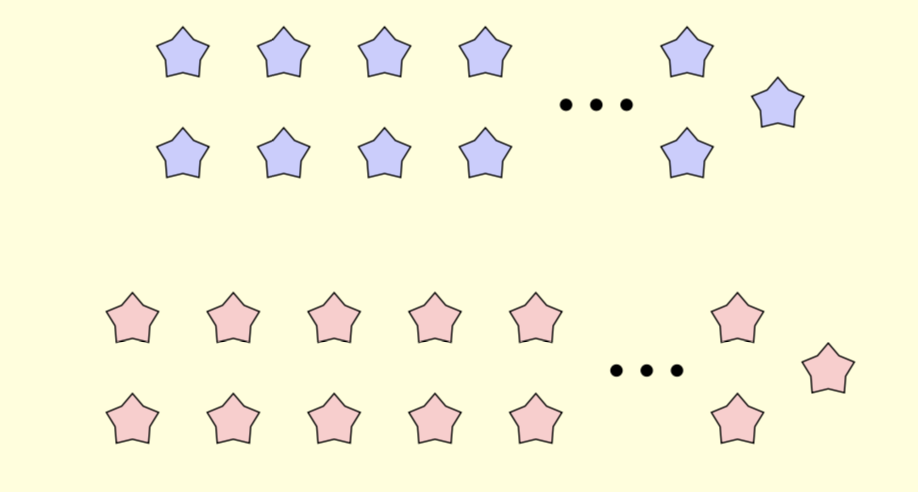}\end{center}

While  such arguments through pictures  lack the formal trappings of proof, 
engaging with such arguments is valuable experience.  Moreover, visual learning of arithmetic 
has a clear basis in the research literature \cite{PaBr13, Ansa16}.  One could place this activity  in a coherent learning
progression across K-12 by revisiting arithmetic of even and odd numbers through
 other arguments, either based on place value and case analysis or grounded in algebra.  Such a learning sequence could advance in
 high-school with activities such as showing that the sum of two squares cannot be one less than a multiple of four.
 
Reasoning about even and odd numbers has been used fruitfully in a number of settings.
Colleagues at the University of Oregon and I use this material as an introduction to different
levels of formalism for undergraduates who aspire to be elementary school teachers \cite{BHS15}.  
Patrick Callahan, a mathematician who led the California Math Project, often has schools evaluate   
student reasoning by asking students at different grade levels to explain why the sum of odd numbers is even.  
Callahan reports that high-schoolers generally fare no better than grade-schoolers, and when presented with the
argument through variables, they (including advanced students)
commonly report that they did not realize it was ``allowed'' for variables to be used in that way.  
 Deborah Ball and her colleagues at the University of Michigan
have used student-driven  discussion about the definition of even numbers as a strong component of teacher training \cite{Bass05, BaBa03}.

The Common Core State Standards for Mathematics were designed through progressions \cite{CCWT18}.
While the Common Core can be read as policy or as informing pedagogy, to a knowledgeable  reader they also suggest proofs for all of K-12
mathematics, short of concepts which require limits, in particular working rigorously with functions over all real numbers.  The commutative
property of multiplication, for example, should be established through noticing that a rectangular array and its transpose are in
bijective correspondence.   Strong curricula engage students in such proofs in age-appropriate ways. 

The Common Core also asks that the canon of elementary mathematics be taught consistently with 
how mathematicians practice mathematics.  The  multiplication table is not just a set of facts but also
a rich locale for conjecture and proof.  Students who type $\frac{1}{7}$ or $\sqrt{5}$ into a calculator and get an ``answer'' 
can  be asked what that answer means, for practice at using definitions as well as reinforcement of estimation and number sense.  
And the story of the law of exponents, which goes from having a 
simple verification for positive whole exponents to being the driver of the definition for all other exponents, is 
a great example of the art of mathematical definition.  (I enjoy the parallel between the law of exponents and the 
homotopy lifting property,
which went from being a property of fiber bundles to being Serre's definition of a fibration.)

Reasoning at the high-school level can reach, for example, the circle of ideas centered on the fact that the sum of the
first $n$ odd numbers is $n^2$.  The statement itself is ripe for conjecture through seeing cases, and it is substantial work
for students to make their conjectures precise.  A graphical argument, as indicated below, is readily accessible, and then 
provides an opportunity to press for details -- for example,
why the number of blocks added in each step increases by two.  Algebraically, there are multiple arguments:  an inductive argument, 
the standard trick for
summing arithmetic sequences, and one
through calculating successive differences (a case in which ``simplifying'' has  a purpose) and employing reasoning central to the
fundamental theorem of calculus.  Indeed, such analysis provides a first explanation, in the discrete setting,
for why quadratic functions should model total displacement in the presence of constant acceleration.

 \begin{center}\includegraphics[keepaspectratio=true, width=0.3\textwidth]{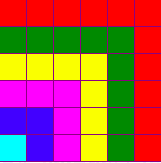}\end{center}

Asking students to provide and critique 
reasoning is more time-intensive than only demonstrating reasoning to them or omitting reasoning altogether.  
But the development of communities of reasoning is sorely needed, especially at this moment.  
Students benefit immensely from establishing truth through logic which 
is accessible to all, with contributions from themselves and peers,
rather than having mathematics exclusively ``handed down'' through authority of teacher and  textbook. 
Such communities foster  responsibility for acknowledging  errors, understanding them as a way to 
progress to correct understanding.

In addition to presenting mathematics  consistently with the practices and values  of our community,
increased reasoning will  pay dividends through deeper retention as well as greater transfer of reasoning ability (see for example
Chapter 3 of \cite{NRC2000}).  
Facility with application is also supported, as students who have access to 
reasoning can more flexibly use mathematics to understand the world.
For example,  concrete visual models are now regularly used to reason about dividing fractions and in particular ``remainders,''
as needed for applications.
If it takes $\frac{1}{4}$ of a ton of steel to make a car, and factory has
$8 \frac{1}{3}$ tons of steel.  Dividing, we see that there could be $35$ cars made,
and  there will be $\frac{1}{12}$ of a ton of steel left over.  This is 
$\frac{1}{3}$ of the amount needed to make another car, as one sees in the quotient $35 \frac{1}{3}$, but in some contexts 
$\frac{1}{12}$ would be the correct ``remainder,'' a quantity which is not accessible if one only learns to ``flip and multiply.''

Progressions occur not only in the mathematics itself, as definitions and theorems build  on previous such, 
but across many domains.  A familiar domain to mathematicians is the amount of abstraction.  Related, but less familiar, is 
the sophistication of representation, including progression from tactile models to pictures to diagrams to variables.  
On the more cognitive side,
there are progressions in student autonomy, including whether tools are called for explicitly
or demanded through problem-solving.  There are progressions in how much of the process of doing mathematics
is engaged in, for example in formulation of conjectures and counterexamples.  Demands in language will also progress.  
Strong curricula attend to all of these types of progression, and address them through research-based pedagogy, 
in particular active learning \cite{PtA14, CBMS16,  MAAIP18}. \\

Teachers need to be prepared to implement curricula which demand that students supply their reasoning.  
As  the mathematical community plays substantial roles in their preparation --
especially that of future high-school teachers -- we 
make some recommendations for teacher training and for other related matters 
below.   For comprehensive recommendations for teacher preparation, see  \cite{METII12, AMTE17}.

\begin{enumerate}

\item Restructure entry-level college courses so students are asked to autonomously provide reasoning,
through proofs or in the context of applications.\footnote{For example, in my department we've replaced our Intermediate Algebra offerings with 
classes which apply algebra to mathematical modeling.}
The rest of the world, including future teachers, takes our choices for these classes as a signal about 
the nature of our subject.  Currently, they deduce that math is only about accurately reproducing procedures!

\item  Develop profound  understanding of the K-12 mathematics progressions in mathematics courses for future elementary and
secondary teachers.

\item Provide math majors interested in teaching with opportunities to connect formal college-level mathematics with school mathematics,
as being developed in multiple NSF-funded projects: the Mathematical Association of America's META Math project,
 the Association of Public and Land-Grant Universities' MODULES project, and the ULTRA project run by Rutgers, Columbia and Temple 
 Universities.
At some point -- perhaps masters-level or other professional preparation-- such students should learn the arguments underlying 
the rules for arithmetic and thus algebra, as preparation for teaching algebra as reasoning and being able to appeal 
to or even fill in background.  Knowledge of such arguments is also helpful for entry-level college teaching.

\item Encourage development of 
opportunities for reasoning in high-school curricula, in particular the use of algebra as a tool of reasoning, for example
to establish divisibility rules or in counting problems.  Reasoning about authentic applications should also be developed, 
and play an especially prominent role in high-school math.

\item Support K-12 educators as they move away  from harmful practices such as acceleration without  
understanding and tracking.  Students should be challenged through depth of understanding, which is achievable in mixed-ability classrooms.
Tracking has  negative consequences for all students, especially those from disadvantaged backgrounds.  
The benefits of such system shifts, as recommended by the National Council of Teachers of Mathematics \cite{NCTM18}, 
 have been borne out by data, for example in the work of San Francisco Unified School District  
\cite{BoFo14, SFUSD18}.  

\end{enumerate}

I would like to thank Yvonne Lai and Jenny Ruef for their help in writing this piece.

\bibliographystyle{alpha}
\bibliography{reasoning}

\newcommand{\etalchar}[1]{$^{#1}$}
\def\cftil#1{\ifmmode\setbox7\hbox{$\accent"5E#1$}\else
  \setbox7\hbox{\accent"5E#1}\penalty 10000\relax\fi\raise 1\ht7
  \hbox{\lower1.15ex\hbox to 1\wd7{\hss\accent"7E\hss}}\penalty 10000
  \hskip-1\wd7\penalty 10000\box7}
  \def\cfudot#1{\ifmmode\setbox7\hbox{$\accent"5E#1$}\else
  \setbox7\hbox{\accent"5E#1}\penalty 10000\relax\fi\raise 1\ht7
  \hbox{\raise.1ex\hbox to 1\wd7{\hss.\hss}}\penalty 10000 \hskip-1\wd7\penalty
  10000\box7}
\begin{thebibliography}{LBHea14}

\bibitem[AMT17]{AMTE17}
{\em Standrds for preparing teachers of mathematics}.
\newblock Available online at {\tt amte.net/standards}. Association of
  Mathematics Teacher Educators, 2017.

\bibitem[Ans16]{Ansa16}
Daniel Ansari.
\newblock The neural roots of mathematical expertise.
\newblock {\em Proceedings of the National Academy of Sciences},
  113(18):4887--4889, 2016.

\bibitem[Bas05]{Bass05}
Hyman Bass.
\newblock Mathematics, mathematicians, and mathematics education.
\newblock {\em Bulletin (New Series) of the American Mathematical Society},
  42(4):417--430, 2005.

\bibitem[BB03]{BaBa03}
Deborah Ball and Hyman Bass.
\newblock Making mathematics reasonable in school.
\newblock In {\em A research companion to principles and standards in school
  mathematics}. National Council of Teachers of Mathematics, 2003.

\bibitem[Bea16]{CBMS16}
Benjamin Braun and et. al.
\newblock Active learning in post-secondary mathematics education.
\newblock Technical report, Conference Board of the Mathematical Sciences,
  2016.

\bibitem[BF14]{BoFo14}
Jo~Boaler and David Foster.
\newblock Raising expectations and achievement. The impact of wide-scale
  mathematics reform giving all students access to high quality mathematics.
\newblock 2014.

\bibitem[BHS]{BHS15}
Tricia Bevans, Andrew Hampton, and Dev Sinha.
\newblock Notes on elementary mathematics.
\newblock Chapter 0 available at
 {\tt https://pages.uoregon.edu/dps/Elementary\_Math\_Chapter\_0.pdf}.

\bibitem[BLS{\etalchar{+}}12]{METII12}
Sybilla Beckmann, W.~James Lewis, Denise Spangler, Alan Tucker, and et~al.
\newblock {\em The mathematical education of teachers, II}, volume~17 of {\em
  CBMS Issues in Mathematics Education}.
\newblock American Mathematical Society, 2012.

\bibitem[{CCSS10}]{CCWT10}
{Common Core Writing Team}.
\newblock Common Core State Standards for Mathematics, available at
{\tt http://www.corestandards.org/Math/}

\bibitem[CCSS18]{CCWT18}
{Common Core Writing Team}.
\newblock Progressions documents, available at {\tt
  http://math.arizona.edu/$\sim$ime/progressions/ }.

\bibitem[NRC00]{NRC2000}
National~Research Council.
\newblock {\em How People Learn: Brain, Mind, Experience, and School: Expanded
  Edition}.
\newblock The National Academies Press, 2000.

\bibitem[KCJS02]{KCSS02}
E.~J. Knuth, Choppin, M.~J., Slaughter, and J.~Sutherland.
\newblock Mapping the conceptual terrain of middle school students'
  competencies in justifying and proving.
\newblock In {\em Proceedings of the 24th annual meeting of the North American
  Chapter of the International Group for the Psychology of Mathematics
  Education}. Clearinghouse for Science, Mathematics, and Environmental
  Education, 2002.

\bibitem[Knu02]{Knu02}
E.~J. Knuth.
\newblock Proof as a tool for learning mathematics.
\newblock {\em Mathematics Teacher}, 95(7):486--491, 2002.

\bibitem[LBHea14]{PtA14}
Steven Leinwand, Daniel Brahier, DeAnn Huinker, and et~al.
\newblock {\em Principles to Actions}.
\newblock National Council of Teachers of Mathematics, 2014.

\bibitem[MAA18]{MAAIP18}
{\em Instructional Practices Guide}.
\newblock Mathematical Association of America, 2018.

\bibitem[NCTM00]{NCTM00}
{\em Principles and Standards for School Mathematics}.
\newblock National Council of Teachers of Mathematics, 2000.

\bibitem[NCTM18]{NCTM18}
{\em Catalyzing Change in High-School Mathematics}.
\newblock National Council of Teachers of Mathematics, 2018.


\bibitem[PB13]{PaBr13}
Joonkoo Park and Elizabeth Brannon.
\newblock Training the approximate number system improves math proficiency.
\newblock {\em Psychological Science}, 24(10):2013--2019, 2013.

\bibitem[SBK09]{SBK10}
D.~Stylianou, M.~Blanton, and E.~Knuth.
\newblock {\em Teaching and learning proof across the grades: A K-16
  perspective.}
\newblock Routledge, 2009.

\bibitem[SFU18]{SFUSD18}
Four years strong in San Francisco: holding an equity-based detracking policy
  over time.
\newblock Technical report, San Francisco Unified School District, 2018.

\end{thebibliography}

\end{document}